\RequirePackage{fix-cm}
\documentclass[twocolumn]{svjour3}					
\smartqed	
\usepackage{graphicx}
\usepackage[T1]{fontenc}
\usepackage[utf8]{inputenc}
\usepackage{listings}
\usepackage{geometry}
\usepackage{fancybox}
\usepackage[cmintegrals]{newtxmath}
\usepackage{bm} 
\usepackage{calc}
\usepackage{amssymb}
\usepackage{stmaryrd}	
\usepackage{color}
\usepackage{amsfonts}
\usepackage{amsmath}
\usepackage[pdfborder={0 0 0},colorlinks=true,linkcolor=blue]{hyperref}
\usepackage{graphicx}
\usepackage{esint}
\usepackage{float}
\usepackage[english]{babel}

\definecolor{commentgris}{rgb}{0.35,0.35,0.35}

\newcommand{\N}{\mathbb{N}}
\newcommand{\Z}{\mathbb{Z}}

\newcommand{\normal}{\mathcal{N}}
\newcommand{\uniform}{\mathcal{U}}
\newcommand{\A}{\mathcal{A}}

\renewcommand{\vec}[1]{\mathbf{#1}}

\newcommand{\type}[1]{{\color{blue}#1}}

\newcommand{\optionCompilation}[1]{{\color{commentgris}#1}}

\newcommand{\prog}[1]{\textit{#1}}

\newcommand{\blue}[1]{{\color{blue}#1}}

\newcommand{\spaceparagraphe}[1]{\vspace{0.5cm}}

\renewcommand{\prod}[1]{\textit{#1}}

\newcommand{\CTA}{CTA}
\newcommand{\panel}[1]{\textbf{#1}}

\begin{document}\sloppy

\title{Polynomial data compression for large-scale physics experiments
}


\author{Pierre Aubert				 \and
				Thomas Vuillaume			\and
				Gilles Maurin				 \and
				Jean Jacquemier			 \and
				Giovanni Lamanna 	\and
				Nahid Emad
}

\authorrunning{P. Aubert, T. Vuillaume, J. Jacquemier, G. Maurin, G. Lamanna \& N. Emad} 

\institute{P. Aubert \and T. Vuillaume \and G. Maurin \and J. Jacquemier \and G. Lamanna \at
							Univ. Grenoble Alpes, Univ. Savoie Mont Blanc, CNRS, LAPP, 74000 Annecy, France \\
							\email{pierre.aubert@lapp.in2p3.fr}					 
					 \and
					 P. Aubert \and N. Emad \at
							Laboratoire d’informatique Parallélisme Réseaux Algorithmes Distribués, UFR des Sciences 45 avenue des États-Unis 78035 Versailles
						\and
					P. Aubert \and N. Emad \at
					Maison de la Simulation, Université de Versailles Saint-Quentin-en-Yvelines, USR 3441 CEA Saclay 91191 Gif-sur-Yvette cedex
}

\date{Received: date / Accepted: date}

\maketitle

\begin{abstract}
The new generation research experiments will introduce huge data surge to a continuously increasing data production by current experiments.
This data surge necessitates efficient compression techniques.
These compression techniques must guarantee an optimum tradeoff between compression rate and the corresponding compression /decompression speed ratio without affecting the data integrity.

This work presents a lossless compression algorithm to compress physics data generated by Astronomy, Astrophysics and Particle Physics experiments.

The developed algorithms have been tuned and tested on a real use case~: the next generation ground-based high-energy gamma ray observatory, Cherenkov Telescope Array (CTA),  requiring important compression performance.
Stand-alone, the proposed compression method is very fast and reasonably efficient.
Alternatively, applied as pre-compression algorithm, it can accelerate common methods like LZMA, keeping close performance.

\end{abstract}

\section{Introduction}
\label{intro}

Several current and next generation experimental infrastructures are concerned by increasing volume of data that they generate and manage. This is also the case in the Astrophysics and Astroparticle Physics research domains where several projects are going to produce a data deluge of the order of several tens of Peta-Bytes (PB) per year \cite{Berghofer:2015glj} (as in the case of CTA) up to some Exa-Bytes (as for the next generation astronomical radio observatory SKA\cite{refSKA}). Such an increasing data-rate implies considerable technical issues at all levels of the data flow, such as data storage, processing, dissemination and preservation.

The most efficient compression algorithms generally used for pictures (JPEG), videos (H264) or music (MP3) files, provide compression ratios greater than $10$.
These algorithms are lossy, therefore not applicable in scientific context where the data content is critical and inexact approximations and/or partial data discarding are not acceptable.
In the context of this work, we focus on compression methods to respond to data size reduction storage, handling, and transmitting issues while not compromising the data content.

Following types of lossless compression methods are applicable for aforementioned situations.
LZMA \cite{refLZMA}, LZ78 \cite{refLZ78}, BZIP2 \cite{refBZIP}, GZIP \cite{refGZIP}, Zstandard \cite{refDataCompressionZstandard} or the Huffman algorithms are often employed because they provide the best compression ratio. The compression speeds of these methods however impose significant constraints considering the data volumes at hand.

Characters lossless compression, CTW (Context Tree Weighting) \cite{refCTW}, LZ77 \cite{refLZ77}, LZW \cite{refLZW}, Burrows-Wheeler transform, or PPM \cite{refPPM}, cannot be used efficiently on physics data as they do not have the same characteristics as text data, like the occurance or repetition of characters.
Other experiments have recently solved this data compression issue \cite{refDataCompressionGAPD}, \cite{refDataCompressionFits} for smaller data rates.

With the increasing data rate, both the compression speed and ratio have to be improved. This paper primarily addresses the data compression challenges.
In this paper, we propose a polynomial approach to compress integer data dominated by a white noise in a shorter time than the classical methods with a reasonable compression ratio.
This paper focuses on both the compression ratio and time because the decompression time is typically shorter.

The paper is organized as follows. Section \ref{secMotivation} explains some motivations.
Section \ref{secPolynomialCompression} describes our three polynomial compression methods.
Section \ref{secImprovingMethod} reports the improvement obtained from our best polynomial compression method on given distributions and CTA data
\cite{articlePresentationCTA}.
Section \ref{secByteOccurrence} gives further details about compression quality.
In section \ref{secConclusion}, some concluding remarks and future plans will be given.

\section{Motivation}
\label{secMotivation}

As the data volumes generated by current and upcoming experiments rapidly increase, the transfer
and storage of data becomes an economical and technical issue. As an example, CTA, the next generation ground-based gamma-ray observatory, will generate hundreds PB of data by 2030. The CTA facility is based on two observing sites, one per hemisphere and will be composed of more than one hundred telescopes in total. Each of them is equipped with photo-sensors equipping the telescopes’ cameras and generating about two hundred PB/year of uncompressed raw data that are then reduced on sites after data selection conditions to the order of the PB/year off-site data yield. The CTA pipeline thus implies a need for both lossy and lossless compression, and the amount of
lossy compression should be minimized while also ensuring good data reading and writing speed. The writing speed needs to be close to real-time, since there is limited capacity on site to buffer such large data volumes. Furthermore, decompression speed is also an issue; the whole cumulated data are expected to be reprocessed yearly, which means that the amount of data needed to be read from disk (, decompressed) and processed will grow each year (e.g. 4 PB, 8 PB, 12 PB, ...).

\begin{figure}[!t]
	\begin{center}
		\includegraphics[width=0.47\textwidth]{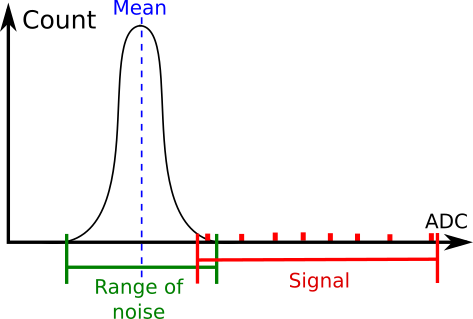}
	\end{center}
	\caption{Example of analog signal digitization in most physics experiments.
		In many cases the white noise (a Gaussian distribution) dominates the signal (generally a Poissonian distribution).
		So, the biggest part of the data we want to compress follows a Gaussian distribution.}
	\label{figBasicWhiteNoiseDistribution}
\end{figure}

In CTA, as in many other experiments, the data acquired by digitization can be described by two components: a Poissonian distribution representing the signal, dominated by a Gaussian-like distribution representing the noise, which is most commonly white noise. As shown in figure \ref{figBasicWhiteNoiseDistribution}, the noise generally significantly dominates the searched signal.

In this paper, we propose a compression algorithm optimised on experimental situation with such characteristics, Gaussian distribution added to a Poissonian one.

Furthermore, in order to respond to time requirement and allow for almost real-time execution the proposed solution can be also combined with the most powerful known compression algorithms such as LZMA to increase tremendously its speed.

\section{The polynomial compression}
\label{secPolynomialCompression}
\begin{figure}[!h]
	\begin{center}
		\includegraphics[width=0.47\textwidth]{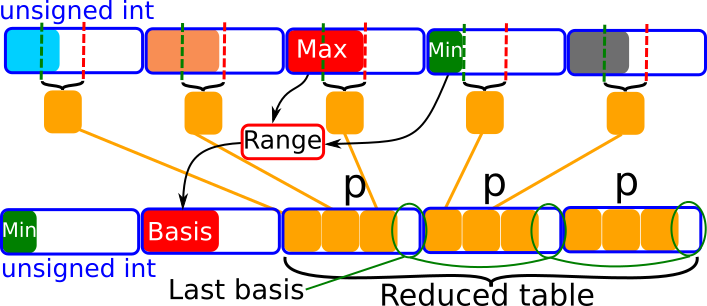}
	\end{center}
	\caption{Illustration of the reduction principle.
		The upper line represents the data (different colours for different values).
		In the second line, the orange blocks represent the changes between the different values to compress.
		The last line shows the compressed data (as they are stored).
		First, the minimum value of the data, next, the base $b = max - min + 1$, which defines the data variations set, $\Z/b\Z$, finally the data variations.
		Several data can be stored in the same \type{unsigned int} and
		only the changes between the data are stored. The common parameters like the range of the data (minimum and maximum or compression base) are stored only once.}
	\label{illustrationSimpleCompression}
\end{figure}

An \type{unsigned int} range, $\left\llbracket 0, 2^{32}\right\llbracket$
defines a mathematical set $\Z/d\Z$, called ring, where $d = 2^{32}$.
The digitized data also define a ring, in this case,
the minimum is $v_{min}$ and the maximum is $v_{max}$ so the corresponding ring is defined as $\Z/b\Z$ with $b = v_{max} - v_{min} + 1$.
In many cases $b < d$, so, it is possible to store several pieces of data in the same \type{unsigned int} (see in figure \ref{illustrationSimpleCompression}).
This compression can be made by using a polynomial approach.
The power of a base is given by the values range.
This allows to add different values in the same integer and compute them back.

\subsection{Basic compression method}
Considering a $n$ elements data vector, $\vec{v} \in \N^N$, its minimum, $v_{min}$ and its maximum $v_{max}$ define its associated ring.
If the data ring is smaller than the \type{unsigned int} ring, it is possible to store several values in one \type{unsigned int}.
The smaller is the base, the higher is the compression ratio.
As the data are in $\left\llbracket v_{min}, v_{max}\right\rrbracket$, the range between $0$ and $v_{min}$ is useless.
Therefore, the data can be compressed by subtracting the minimum value, forming a smaller base.
The minimum can be stored once before the compressed data.
The compression base $B$ is defined by : $B = v_{max} - v_{min} + 1$.
With this base we are able to store $\left(v_{max} - v_{min}\right)$ different values.
The compression ratio, $p$, is given by the number of bases $B$ that can be stored in one \type{unsigned int}
(in $\left\llbracket 0, 2^{32} \right\llbracket $)~:
\begin{eqnarray}
	p & = & \left\lfloor \dfrac{\ln \left( 2^{32} - 1 \right)}{\ln B} \right\rfloor \label{firstEqP}
\end{eqnarray}

The compressed elements, $s_j$, are given by~:

\begin{eqnarray}
	s_j & = & \sum_{i=1}^p v_{i + p\times (j - 1)} \times B^{i-1} \quad \text{for } 1 \leq j < \frac{n}{p}
\end{eqnarray}

A polynomial division can be used to uncompress the data.

\subsection{Advanced compression method}
\begin{figure}
	\begin{center}
		\includegraphics[width=0.47\textwidth]{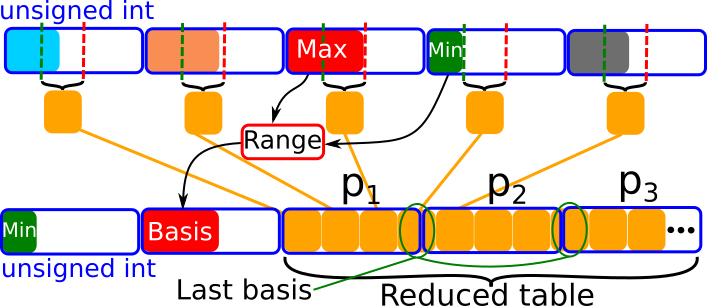}
	\end{center}
	\caption{Illustration of the advanced reduction.
	The upper line represents the data (different colours for different values).
	In the second line, the orange blocks represent the changes between the different values to compress.
	The last line shows the compressed data (as they are stored).
	First, the minimum value of the data, next, the base $b = max - min + 1$, which defines the data variations set, $\Z/b\Z$, and finally the data variations.
	The storage space is optimized by avoiding useless gaps between data.
	With this method there is no useless space to store the compressed data.}
	\label{illustrationCompressionBasis}
\end{figure}

The inconvenience of the basic compression method is the unused space at the end of each packed \type{unsigned int} (see in figure \ref{illustrationSimpleCompression}).
The ideal case is the one that has no unused space when storing the compressed
data illustrated in figure \ref{illustrationCompressionBasis} which avoid unused space.
The advanced polynomial compression tends to become ideal case.
However, minimizing the time for read and write provides a faster compression and decompression speed.

\begin{figure}
	\begin{center}
		\includegraphics[width=0.47\textwidth]{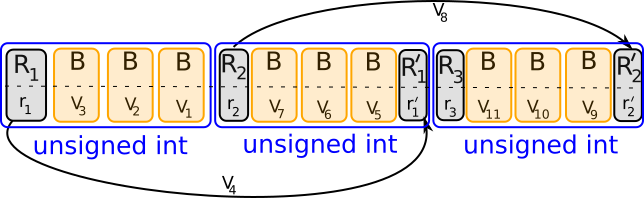}
	\end{center}
	\caption{This figure shows how the data of the vector $\vec{v}$ are stored in the packed vector.
	The first line gives the base used to store the values, the second line shows the variables used to store the values with respect to their base.
	To increase the compression ratio we need to split the last base $B$ in to base $R$ and $R^\prime$ in order to use the storage capacity of an \type{unsigned int} as much as we can.
	The values on the left of an \type{unsigned int} are stored with a low power of the base $B$.
	The values on the right of an \type{unsigned int} are stored with a high power of the base $B$.}
	\label{figSplittingBaseB}
\end{figure}

Compression ratio can be improved by splitting the last base (see figure \ref{illustrationSimpleCompression} and figure \ref{illustrationCompressionBasis}),
used to pack less data in the
same \type{unsigned int}, into two other bases, $R$ and $R^\prime$ in order to have $B \leq R\times R^\prime$.
In this case, the base $R$ is stored in the current packed \type{unsigned int} and the base $R^\prime$ is stored in the next one (see in figure \ref{figSplittingBaseB}).
This configuration ensures a more efficient data order for CPU data pre-fetching at the decompression time, in order to ensure decompression faster than compression.

The data are accumulated from the highest exponent of the base $B$ to the lowest.
This ensures the decompression will produce uncompressed contiguous data.

This splitting stores a value to be compressed, $v$, in two bases, $R$ and $R^\prime$ with two variables $r$ and $r^\prime$.
The variables $r$ and $r^\prime$ are stored in two consecutive packed element (\type{unsigned int}).

The calculation of the bases $R$ and $R^\prime$ is possible when the number of bases $B$ that can be stored in an \type{unsigned int} is known.
The number of bases $B$ that can be stored in the first packed \type{unsigned int}, $p_1$, is given by the following equation :

\begin{eqnarray}
	p_1 & = & \left\lfloor \dfrac{\ln \left( 2^{32} - 1 \right)}{\ln B} \right\rfloor
\end{eqnarray}

The split base $R_1$ is given by :

\begin{eqnarray}
	R_1 & = & \left\lfloor \dfrac{2^{32} - 1}{B^{p_1}} \right\rfloor
\end{eqnarray}

The $R_1$ base must be completed, to store an element $e \in \left\llbracket 0, B \right\llbracket $, by :

\begin{eqnarray}
	R^\prime_1 & = & \left\lfloor \dfrac{B}{R_1} \right\rfloor + \left( 1 \text{ if } B\mod R_0 \neq 0 \right) \label{baseRPrimeCalculation}
\end{eqnarray}

\newpage{}

So $R_1 \times R^\prime_1 \geq B$. Each base $R_i$ and $R^\prime_i$ are associated to a stored value $r_i$ and $r^\prime_i$ respectively.
The first packed element $s_1$ can be written as follow :

\begin{eqnarray}
	s_1 & = & r_1 + R_1\times \left( \sum_{k=1, k \neq p_1}^{p_1} v_k B^{p_1 - k} \right)
\end{eqnarray}

Where :

\begin{eqnarray}
	r^\prime_1 & = & \left\lfloor \dfrac{v_{p_1}}{R_1} \right\rfloor \\
	r_1	 & = & v_{p_1} - r^\prime_1 \times R_1
\end{eqnarray}

The value $r_1$ is associated to the base $R$ and the value $r^\prime_1$ is associated to the base $R^\prime$.

The number of bases $B$ that can be stored in the second packed element, $p_2$, is given by :

\begin{eqnarray}
	p_2 & = & \left\lfloor \dfrac{\ln \left( 2^{32} - 1 \right) - ln R^\prime_1 }{\ln B} \right\rfloor
\end{eqnarray}

The rest split base $R_2$ can be written as :

\begin{eqnarray}
	R_2 & = & \left\lfloor \dfrac{2^{32} - 1}{R^\prime_1 B^{p_2}} \right\rfloor
\end{eqnarray}

The base $R_2$ the second packed element can be calculated~:
\begin{eqnarray}
	s_2 & = & r_2 \nonumber\\ & & + R_2\times \left( \sum_{k=1, k \neq p_1}^{p_1} v_{p_0 + k} B^{p_1 - k} + r^\prime_1	B^{p_1} \right)
\end{eqnarray}

The equation \ref{baseRPrimeCalculation} can be used to calculate $R^\prime_2$.\\

The compression of a whole vector can be done by using a mathematical series to calculate the split base for each packed element.
Assuming the base $R^\prime_0 = 1$ for the first step, the mathematical series used to compress an entire vector of \type{unsigned int} can be written as follow (for $0 < i \leq n_p$, where $n_p$ is the number of packed elements) :

\begin{eqnarray*}
	p_i		& = & \left\lfloor \dfrac{\ln \left( 2^{32} - 1 \right) - ln R^\prime_{i-1} }{\ln B} \right\rfloor		\\
	R_i		& = & \left\lfloor \dfrac{2^{32} - 1}{R^\prime_{i-1} B^{p_i}} \right\rfloor			\\
	R^\prime_i	& = & \left\lceil \dfrac{B}{R_{i}} \right\rceil						\\
	q_i		& = & i - 1 + \sum_{k=1, k \neq i}^i p_k	\quad \text{ or } 0 \text{ if } i = 0	\\
	r^\prime_i & = & \left\lfloor \dfrac{v_{q_i + p_i}}{R_i} \right\rfloor \\
	r_i	 & = & v_{q_i + p_i} - r^\prime_i \times R_i			\\
	s_i		& = & r_i + R_i\times \left( \sum_{k=1, k \neq p_i}^{p_i} v_{q_i + k + 1} B^{p_i - k} + r^\prime_{i-1}	B^{p_i} \right)
\end{eqnarray*}

Where $p_i$ is the number of bases $B$ that can be stored in the $\text{i}^\text{th}$ packed element, $R_i$ and $R^\prime_i$ are the split base, $r_i$ and $r^\prime$ their corresponding values,
$q_i$ is used to know how many elements have been packed until the $\text{i}^\text{th}$ packed element, finally $s_i$ is the value of the $\text{i}^\text{th}$ packed element.

\subsection{Blocked compression method}

We observed that smaller the signal range, more efficient is the compression.
The advanced compression method presented above is particularly efficient to compress white noise with small spread.
Conversely, if the gaussian noise or the poissonian signal is spread out, the efficiency decreases.
However, the efficiency can be improved by dividing the vector into less items to diminish the impact by the higher values on the global compression ratio.

The block efficiency, and their size determination will be discussed in the section \ref{subSecReductionGivenDistri}.

\section{Experiments and analysis}\label{secImprovingMethod}

\begin{figure}[!t]
	\begin{center}
		\includegraphics[width=0.46\textwidth]{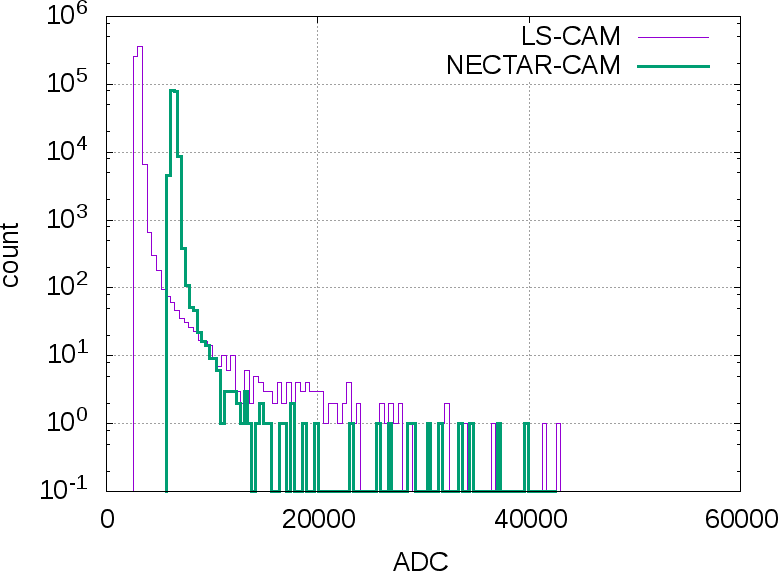}
	\end{center}
	\caption{This figure illustrates the typical signal distribution obtained in several of the cameras used in \CTA{} \cite{articlePresentationCTA}.
	}
	\label{figDistributionCtaCam}
\end{figure}

In order to test and evaluate the performance of the previously
described polynomial compression method, in the following Monte Carlo
simulated distributions will be used. These distributions are in
agreement with the measured data from Cherenkov cameras (see figure \ref{figDistributionCtaCam}).

	\subsection{Simulation of the distribution}

		The data can be described by a random gaussian distribution with a given standard deviation
		(the white noise in the cameras' signals)
		and by adding a uniform distribution in the given camera's signal range (the physics signal).
		Consider the set of the camera pixels distribution values, $\A$~:
		\begin{equation}
			\A\left(\mu, \sigma, x, y, a, N\right) = \normal\left(\mu, \sigma \right)^{N-a} \cup \uniform\left(x, y\right)^a
		\end{equation}

		Where :
		\begin{itemize}
			\item $\mu$ : gaussian noise mean
			\item $\sigma$ : gaussian noise standard deviation
			\item $\left(x, y\right)$ : range of uniform signal value
			\item $a$ : number of values in the uniform distribution (signal)
			\item $N$ : total number of values in the vector
		\end{itemize}
		($\normal$ is the normal distribution, the simulated noise, and the $\uniform$ describes an uniform distribution, the simulated signal.)

\begin{figure}[!h]
	\begin{center}
		\includegraphics[width=0.47\textwidth]{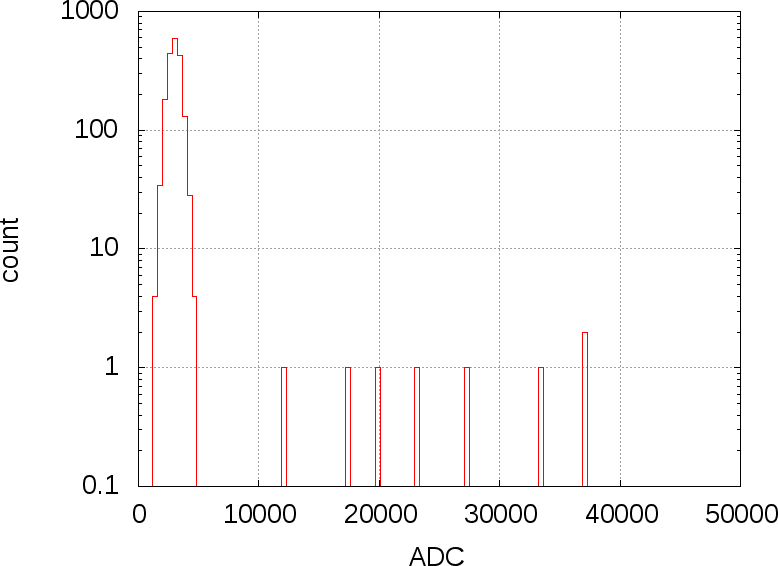}
	\end{center}
	\vspace{-0.4cm}
	\caption{Typical simulated distribution used to improve the data reduction in the set $\A\left(3\,000, 500, 2\,000, 45\,000, 9, 1855\right)$.
		In this case, the gaussian distribution represents $99\%$ of the pixels' values and the uniform distribution represents $1\%$ of this distribution.}
	\label{simulatedDistributionCameraAdc}
\end{figure}

		An example of simulation is presented on figure \ref{simulatedDistributionCameraAdc}.

		In this paper we have tested the distribution $\A$ with $\mu = 3000$, $\sigma \in \left[100, 10\,000\right]$,
		$x = 2\,000$, $y \in \left\llbracket 20\,000, 100\,000\right\rrbracket$,
		$a \in \left\llbracket 1, 1\,000\right\rrbracket$ and $N \in \left\llbracket 1855, 10\,000\right\rrbracket$.

\vspace{-0.5cm}
	\subsection{Polynomial reduction on given distributions}\label{subSecReductionGivenDistri}

		The implementation of blocked polynomial reduction has been tested on given distributions.
		This test determines the influence of the distribution parameters on the compression ratio.

\begin{figure}[!h]
\begin{center}
	\includegraphics[width=0.47\textwidth]{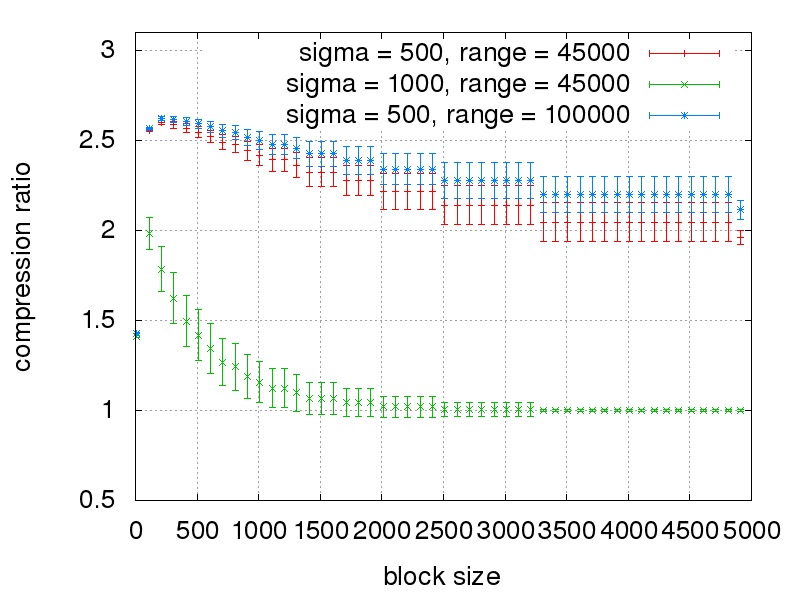}
	\includegraphics[width=0.47\textwidth]{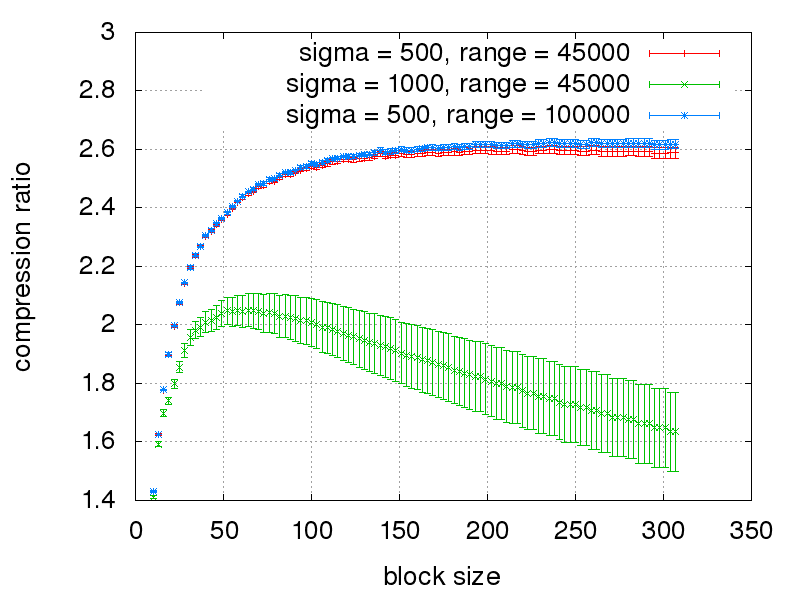}
\end{center}
\vspace{-0.5cm}
\caption{\panel{Top panel} : illustration of the compression ratio versus the number of element in the blocks used to compressed a vector of data.
	The red points (+) give the compression ratio for a distribution with $\sigma = 500$ and $range = 45\,000$,
	in $\A\left(3\,000, 500, 2\,000, 45\,000, 4, 10\,000\right)$.
	The blue points ($\ast$) give the compression ratio for a distribution with $\sigma = 500$ and $range = 100\,000$,
	in $\A\left(3\,000, 500, 2\,000, 100\,000, 4, 10\,000\right)$.
	The green points ($\times$) give the compression ratio for a distribution with $\sigma = 1\,000$ and $range = 45\,000$,
	in $\A\left(3\,000, 1\,000, 2\,000, 45\,000, 4, 10\,000\right)$.
	The tails of the plots give the compression ratio of the advanced polynomial reduction method.
	\panel{Bottom panel} : the same plot zoomed.}
\label{cmpErrorSimulatedDist4}
\end{figure}

		As the polynomial compression uses statistical properties to compress data, the test can only be done with a set of distributions.
		The figure \ref{cmpErrorSimulatedDist4} shows the compression ratio for $1000$ vectors with
		$\A\left(3000, 500, 2000, 45\,000, 4, 1855\right)$ to compute the variations (red curve).
		The gaussian $\sigma$ variation has a high influence on the final compression ratio, of the order of $25\%$
		from $\sigma = 1000$ to $\sigma = 500$ in the best block size case.
		The signal range influence is lighter, $5\%$ or $10\%$ depending on the block size, and $5\%$ for the best block size.
		The block size choice is important too. The compression ratio is weaker if the blocks are too long,
		$25\%$ of lower compression for $\sigma = 500$ and $30\%$ for $\sigma = 1000$.
		In this case using blocks of $154$ elements allows a compression ratio of $2.47054$ which is larger than $17\%$ to the basic compression ($2.10361$).

\vspace{-1cm}

\begin{figure}
	\begin{center}
		\includegraphics[width=0.47\textwidth]{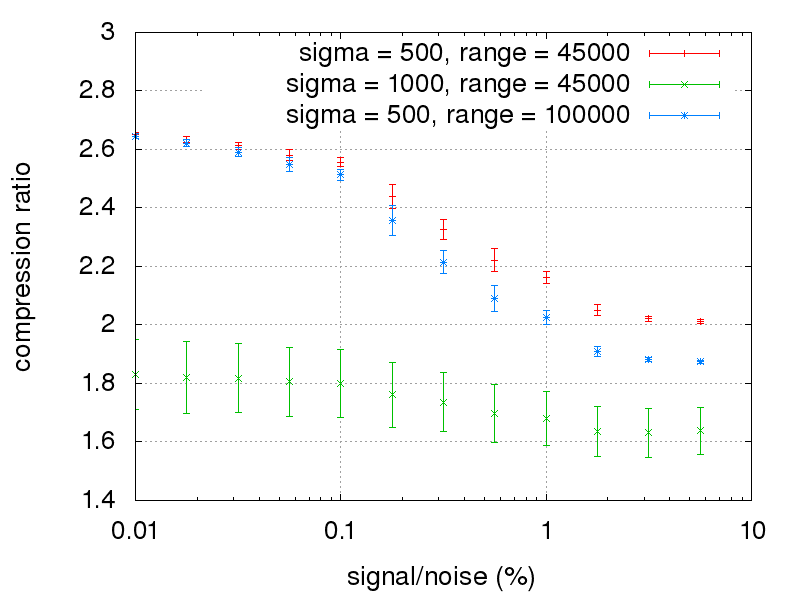}
	\end{center}
	\vspace{-0.5cm}
	\caption{Illustration of the signal/noise ratio influence on the final compression ratio for a vector of $10\,000$ values.
		The red points (+) give the compression ratio for a distribution with $\sigma = 500$ and $range = 45\,000$.
		The blue points ($\ast$) give the compression ratio for a distribution with $\sigma = 500$ and $range = 100\,000$.
		The green points ($\times$) give the compression ratio for a distribution with $\sigma = 1\,000$ and $range = 45\,000$.
		The compression ratio is constant from $10\%$.}
	\label{figCmpRatioNbUniformOnNbPixel}
\end{figure}

\vspace{-0.5cm}
Figure \ref{figCmpRatioNbUniformOnNbPixel} shows that the signal/noise ratio has high influence
on the compression ratio if the spread is not too high.
On the contrary, if the noise spread is important (green curve), the signal/noise ratio has less influence on the final compression ratio.

	\vspace{-0.5cm}
	\subsection{Polynomial reduction on CTA data}

		In the previous section \ref{subSecReductionGivenDistri} we have described the compression ratios obtained with the blocked polynomial reduction applied on modelled/simulated data distributions.
		Such an improvement cannot reflect properly the compression ratio with physics data, since they result typically in a superposition of several distributions coming from several photo-sensors (e.g. pixels) read simultaneously. 
	
\begin{table*}[!t]
\caption{The polynomial compression ratio, time and compressed file size compare the LZMA (best compression existing).
The tested file is the full waveform simulation of the PROD\_3 (run 3998) of the CTA experiment.
The used CPU was a Intel core i7 M 560 with 19 GB of RAM installed with a Fedora .
}
\label{ratioTimeCompressionPolynomialSliceLZMA}
\centering
\scalebox{0.68}{
\begin{tabular}{|c|c|c|c|c|c|c|}
	\hline
					&	Compression	&	Compression			&	File size (GB)	& Decompression		&	Compression	 &	Decompression	\\
					&	ratio		&	Elapsed Time			&			& Elapsed Time		&	Elapsed Time RAM &	Elapsed Time RAM	\\
	\hline
	\hline
	No compression			&	$1$		&	$0$				&	$7.6$		&	$0$		&	$0$		&	$0$	\\
	\hline
	\hline
	Advanced Polynomial		&	$2.71$		&	$5\,$min$36.025$s		&	$2.8$		&	$2\,$min$35\,$s	&	$2\,$min$56\,$s	&	$2\,$min$30\,$s	\\
	Reduction			&			&					&			&			&	&	\\
	\hline
	BZIP2				&	$2.62$ 		&	$19\,$min$\,18.247\,$s		&	$2.9$		&	$2\,$min$23.676\,$s	&	- 	&	-\\
	\hline
	LZMA (-mx=9 -mfb=64 -md=200m)	&	$6.52$		&	$2\,$h$14\,$min$44\,$s		&	$1.166$		&	$1\,$min$49.689\,$s	&	$2\,$h$03\,$min$44\,$s & $2\,$min$03\,$s\\
	LZMA (-mx=1 -mfb=64 -md=32m)	&	$5.88$		&	$14\,$min$00\,$s		&	$1.293$		&	$2\,$min$10\,$s 	&	$15\,$min$42\,$s	&	$2\,$min$09\,$s	\\
	\hline
	Poly + LZMA (9 16 32)		&	$5.02$		&	$10\,$m$49\,$s			&	$1.513$		&	$2\,$min$13\,$s	&	$4\,$min$57\,$s	&	$2\,$min$15\,$s\\
	\hline
\end{tabular}
}
\end{table*}
	
		We have therefore tested our compression method on Monte Carlo simulated CTA-like data (i.e. Cherenkov light emitted by atmospheric electromagnetic showers and captured by cameras on telescopes).
		Only shower pictures registered in stereoscopy by more than one telescope are recorded.
		Each telescope’s camera produces individually a file containing its own picture (called also an “event”) resulting from the different signals registered by all pixels/photo-sensors of the camera itself.
		An event data-file is then composed of a header, used to describe properties like its timestamp, plus the recorded camera data, e.g. either the integrated signal from all pixels and/or the dynamical evolution of the signal in time (waveform).
		The event data file can have a size of typically several tens of thousands bytes.
		A selection on the pixels to be saved will likely be applied in the acquisition pipeline in order to reduce the final data rate.

		Among the various specifications that have to be fulfilled by the CTA data format, each data file has to be readable by part to enable the access to its header without requiring a full decompression step.
		Therefore, a blocked compression is allowed.
	
\begin{figure}[H]
	\begin{center}
		\includegraphics[width=0.47\textwidth]{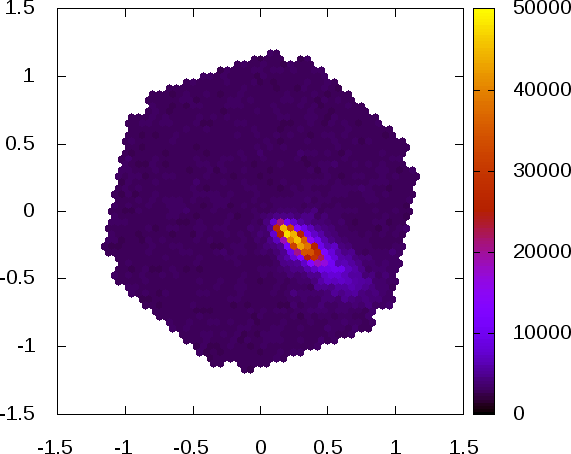}
	\end{center}
	\caption{Illustration of the ellipse shape of a particles shower recorded by a camera in the CTA Monte-Carlo. The color scale represents the number of photons detected in a pixel.}
	\label{lstTelescopeImage}
\end{figure}

		Our test was performed on CTA Monte-Carlo Prod 3 \cite{articleMonteCarlo} files, which simulate telescopes observing the Cherenkov light emitted by particles’ showers in the atmosphere and used to characterize the scientific output of CTA in experimental conditions, thus they are reasonably realistic.
		The CTA Monte Carlo data are converted into a specific high performance data format, which stores pixels’ values in $16\,$ bits to enable fast computation.
		The test files contain $624$ telescopes, $22\,000\,$images, $7.6\,$GB data in waveform mode and $474\,$MB of integrated signal.
		Each image is composed of the lighted pixels concerned by both the noise as well as by the genuine signal and having an almost elliptical shape (see figure \ref{lstTelescopeImage}).
	
		In the following we present the way to adapt the polynomial compression method to the CTA prerequisites.
		All tests are executed on an Intel core i5 clocked at 2.67 GHz with SSE4 instructions without SSD disk.
		

		\vspace{-0.3cm}
		\subsubsection{Test on waveform \CTA{} data}\label{subSubSecPolyWave}
			
			The waveform-data of CTA record the electromagnetic showers’ expansion. Thus, each pixel has values in time.
			The number of values depends on the camera type.
			Each value is digitized in $12$ or $16$ bits and is stored in $16$ bits for computing reasons.
			
		\vspace{-0.7cm}
			To improve waveform compression and enable High Performance Computing (CPU data pre-fetching and vectorization) we choose to store them into matrices. 
			The matrix element $M_{i,j}$ corresponds to the i$^\text{th}$ time of the j$^\text{th}$ pixel of the current camera.
			The row alignment enables a better compression because it increases the number of sequences of similar values.
			This configuration enables the optimisation of waveforms' pictures computing. 
			
			Our polynomial compression reduces the test file of $7.6\,$GB in a file of $2.8\,$GB (compression ratio of $2.71$) in $5\,$min$\,36\,$s.
			The table \ref{ratioTimeCompressionPolynomialSliceLZMA} compares our results with classical methods.
			We test BZIP2 and LZMA algorithms over the whole file.
			
\begin{table*}[!h]
\caption{The polynomial compression ratio, time and compressed file size compare the LZMA (best compression existing).
The tested file is the simulation the PROD\_3 (run 497) of the CTA experiment.
The combination of our advanced polynomial compression and LZMA allows a compression as good as a pure LZMA compression but $19$ times faster.
The used CPU was a Intel core i5 M 560 with 8 GB of RAM installed with an Ubuntu 16.4.
}
\label{ratioTimeCompressionPolynomialLZMA}
\centering
\scalebox{0.83}{
\begin{tabular}{|c|c|c|c|c|c|}
	\hline
					&Compression	&	Elapsed			&	File size (MB)	&	Compression	 &	Decompression	\\
					&ratio		&	Time			&			&	Elapsed Time RAM &	Elapsed Time RAM	\\
					&		&				&			&			&		\\
	\hline
	\hline
	No compression			&	$1$	&	$0$			&	$474$		& $0$ & $0$\\
	\hline
	\hline
	Advanced Polynomial		&	$3.74$	&	$3.7\,$s		&	$127$		& $0.9\,$s & $0.9\,$s\\
	Reduction			&		&				&			& & \\
	\hline
	BZIP2				&	$4.69$ 	&	$1\,$min$\,48\,$s	&	$101$		& $1\,$min$0\,$s & $6.23\,$s\\
	LZMA (7z)			&	$4.84$	&	$7\,$min$\,48.636\,$s	&	$98$		& $1\,$min$\,18\,$s & $9.28\,$s\\
%
	\hline
	\textbf{Advanced Polynomial}	&	\blue{$\mathbf{4.84}$}	&	$\mathbf{24.646\,}$\textbf{s}	&	$\mathbf{98}$	& $1\,$min$\,20\,$s & $11\,$s\\
	\textbf{Reduction + LZMA}	&				&					&			& & \\
	\hline
\end{tabular}
}
\end{table*}

			Our method is $3.5$ times faster than the BZIP2 algorithm and offers a better compression.
			
			For the LZMA algorithm, the program \prog{7z} is used on two tests.
			
\vspace{-0.6cm}
\begin{figure}[!h]
	\begin{center}
		\includegraphics[width=0.47\textwidth]{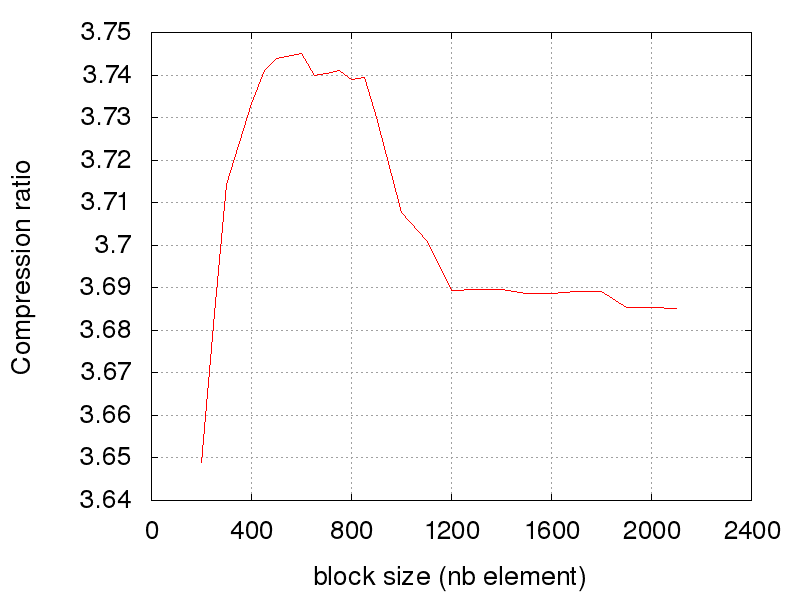}
	\end{center}
	\vspace{-0.5cm}
	\caption{Comparison of the compression ratio for different compression block sizes on the CTA PROD\_3 Monte-Carlo.}
	\label{compresionRatioRun497VsBlockSize}
\end{figure}
		\vspace{-0.6cm}
			
			First we investigate to reach the best compression ratio with the commande line \optionCompilation{7z a -t7z -m0=lzma -mx=9 -mfb=64 -md=200m -ms=on}.
			Thus, the test file is compressed with a compression ratio of $6.52$ in $2\,$h$14\,$min $44\,$s.
			We also investigate the fast compression mode of LZMA (\optionCompilation{7z a -t7z -m0=lzma -mx=1 -mfb=64 -md=32m -ms=on}) and we obtain a compression ratio of $5.88$ in $14\,$m $00\,$s.
			Finally we combine the polynomial compression with LZMA and obtain a compression ratio of $5.02$ in only $10\,$m $49\,$s.
			The decompression times of the different algorithms are similar.
			However, compression/decompression over a whole file are not in agreement with the data requirement (events or blocks have to be compressed separatelly).
			A solution is to compress several events packed in blocks with a higher granularity level.
			With LZMA (native), this method archieves a compression ratio of $6$ by compressing $300$ events per block which represents less than $0.02$ second of signal for the LST-CAM.
			In this case, each block contains approximately $100\,$MB of data (depending on the camera).
			
		\subsubsection{Test on integrated CTA data}

	\vspace{-0.2cm}
	
	The integrated data are obtained by the reduction of the waveform signal.
	This reduction can be performed on all the pixels' waveforms or on several pixels' waveforms.
	In our case, we reduced matrices of section \ref{subSubSecPolyWave} in vectors to enable High Performance Computing.
	The produced vector can be described by $24$ or $32$ bits data depending on the cameras.
	
	The polynomial compression archieves a compression ratio of $3.74$ in $3.7\,$s on a $474\,$MB test file.
	
	The figure \ref{compresionRatioRun497VsBlockSize} shows the compression ratios obtained with different compression block size.
	The plot variations denote the different compression ratios from the different cameras in the file.
	Statistically, the camera data have not the same compression ratio.
	This is why there are fluctuations.

	The difference with the simulated distribution comes from the $7$ types of cameras.
	Each camera has a typical ellipse size.
	If one block contains the full ellipse signal, it is less compressed compared to others.
	The result is a better global compression.
	
	Table \ref{ratioTimeCompressionPolynomialLZMA} compares polynomial compression with LZMA and BZIP2 algorithms.

	The BZIP2 algorithm provides a better compression ratio ($4.69$) but in $1\,$min$48\,$s ($29$ times slower than our compression).
	The best compression ratio is obtained with the LZMA algorithm ($4.84$) but in $7\,$min$48\,$s ($126$ times slower than our compression).
	
	By combining our advanced polynomial reduction with the LZMA compression we obtain the same compression ratio as pure LZMA, but $19$ times faster.
	Moreover, the use of the polynomial reduction allows the LZMA to keep its flat profile because it packs small values in high values.
	So the values average increases and the bytes profile becomes more flat.
	
	The polynomial reduction allows also a better compression ratio than a classical bit-shifted compression because the left space in an \type{unsigned int} is used.
	
	Extrapolating to CTA yearly data rate of $4\,$PB \cite{refCtaTechnicalDesignReport}, the usage of the LZMA algorithm in this case would require more than 1750 core.year only for the compression.

	\spaceparagraphe

	\section{Bytes occurrences}\label{secByteOccurrence}
	
	For further test purpose one can compare the distributions of the byte-values as in the initial file and in the compressed one with the purpose of verifying that the compression algorithm has not altered the physical distributions.
	At this point the file cannot be further compressed with a lossless compression \cite{refLZMA}.
	
\begin{figure}
	\begin{center}
		\includegraphics[width=0.47\textwidth]{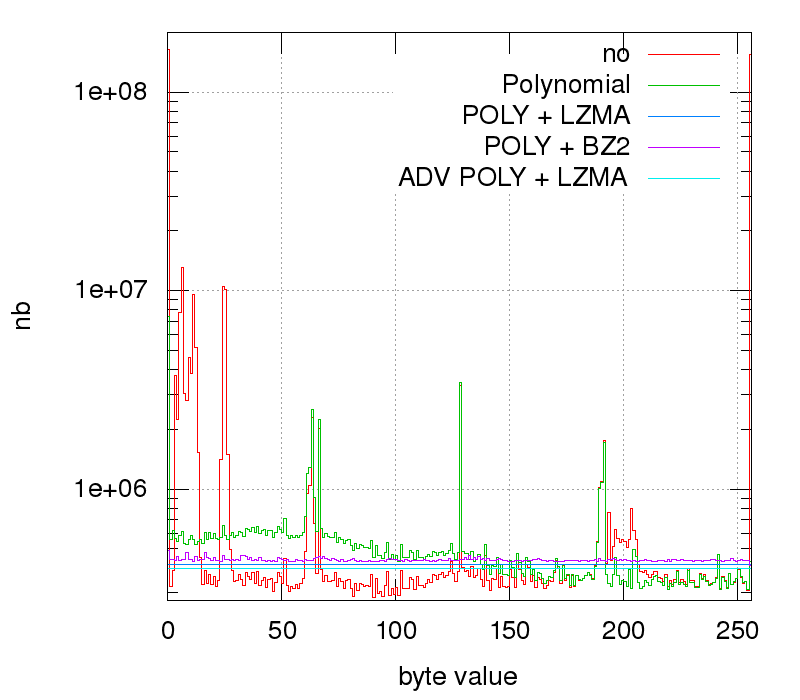}
	\end{center}
	\caption{Comparison of the different values of bytes in a binary file.
		In red, the profile of the uncompressed file (CTA PROD\_3 Monte-Carlo)
		In green, the profile of the only polynomial reduced file.
		In blue, the profile of the corresponding compressed file with the polynomial reduction and LZMA compression.
		In purple, the profile of the corresponding compressed file with the polynomial reduction and GZIP compression.
		In cyan, the profile of the corresponding compressed file with the advanced polynomial reduction and LZMA compression.}
	\label{profileProdFilePolyReductionPolyLZMAPolyGZ}
\end{figure}

	Figure \ref{profileProdFilePolyReductionPolyLZMAPolyGZ} shows the different byte-value distributions in the initial file and in the file compressed with a polynomial reduction, LZMA and BZIP2 for integrated test files.
	
	This figure shows the polynomial compression smoothens the profile.
	The best compression is obtained with the LZMA compression on a file compressed with a polynomial reduction because its profile is flat.
	The combination of the polynomial compression and the BZ2 algorithm does not provide a better compression ratio or a faster compression speed.

\section{Conclusion}\label{secConclusion}

	In this article, we introduced a new lossless compression algorithm to compress integers from digitized signals and dominated by a white noise.

	This method is very fast, helps CPU data pre-fetching and eases vectorization.
	It can be integrated in each data format that deals with tables or matrices of integers.
	This method compresses preferentially integers but an adaptation of this algorithm can enable floating-point data compression with fixed precision that returns integers.
	It compresses matrices and tables separately in order to keep a similar data structure between compressed and uncompressed data.
	The decompression is roughly twice faster than the compression.
	This method can also be vectorized to improve its speed.
	
	Tests on CTA Monte-Carlo data show that the polynomial compression is less efficient than LZMA but more efficient than BZIP2 on waveforms data.
	
	The integrated data compression is very efficient and fast.
	Used as a pre-compression for LZMA, we obtain the same compression ratio as pure LZMA but in a compression duration $19$ times shorter.
	
	The method's simplicity offers easy development in many languages and the possiblity to be used on simple embedded systems,
	or to reduce the data volume produced by high sensitive captors on FPGA.
	It can be also used as pre-compression of stronger methods (like LZMA) and accelerate it.


\begin{acknowledgements}

This work is realised under the Astronomy ESFRI and Research Infrastructure Cluster (ASTERICS project)
supported by the European Commission Framework Programme Horizon 2020
Research and Innovation action under grant agreement n. 653477.
\end{acknowledgements}


The source code of the polynomial compression method discussed in this work is avaliable under \small{\url{https://gitlab.in2p3.fr/CTA-LAPP/PLIBS_8}}.

\end{document}